\documentclass[a4paper,11pt]{article}
\usepackage{latexsym,amssymb,amsmath,textcomp}
\usepackage[english]{babel}

\usepackage{amsthm}
\theoremstyle{plain}
\usepackage{mathtools}

\newtheorem{Theorem*}{Theorem}

\newtheorem{Lemma}{Lemma}
\newtheorem{Definition}{Definition}
\newtheorem{Corollary}{Corollary}
\newtheorem{Proposition}{Proposition}

\usepackage{hyperref}
\usepackage{tikz}
\DeclareMathOperator{\rank}{srank}

\title{A gap in the slice rank of $k$-tensors}
\author{Simone Costa\thanks{DII/DICATAM - Sez. Matematica, Universit\`a degli Studi di Brescia, Via
Branze 38, I-25123 Brescia, Italy. email: simone.costa@unibs.it}\ \
and Marco Dalai%
\thanks{DII, Universit\`a degli Studi di Brescia, Via
Branze 38, I-25123 Brescia, Italy. email: marco.dalai@unibs.it}}\begin{document}
\maketitle
\begin{abstract}
The slice-rank method, introduced by Tao as a symmetrized version of the polynomial method of Croot, Lev and Pach  and Ellenberg and Gijswijt, has proved to be a useful tool in a variety of combinatorial problems. Explicit tensors have been introduced in different contexts but little is known about the limitations of the method. 

In this paper, building upon a method presented by Tao and Sawin, it is proved that the asymptotic slice rank of any $k$-tensor in any field is either $1$ or at least $k/(k-1)^{(k-1)/k}$. This provides evidence that straight-forward application of the method cannot give useful results in certain problems for which non-trivial exponential bounds are already known. An example, actually a motivation for starting this work, is the problem of bounding the size of trifferent sets of sequences, which constitutes a long-standing open problem in information theory and in theoretical computer science. 
\end{abstract}
\section{Introduction}

The polynomial method has been receiving renowed attention since the breakthrough result of Croot, Lev and Pach \cite{croot-lev-pach-2016} and subsequent follow-up results, among which the notable ones by Ellenberg and Gijswijt \cite{ellenberg-gjswijt-2017} and Naslund and Sawin \cite{naslund-sawin-2017}. A useful symmetrized formulation of this particular application of the polynomial method was provided by Tao in \cite{Blog1} based on a notion of \emph{slice rank} of tensors. In this formulation, the size of combinatorial structures under study is upper-bounded by the slice rank of appropriately constructed tensor powers. The notion of slice rank can be interpreted in a more general framework of tensor-ranks which is given an in-depth discussion in \cite{jeroen-2018}.

The slice rank method has been applied to several combinatorial problems such as the tri-colored sum-free sets, the sunflowers free sets, the capsets and the progression-free problem. In those cases, the method gave the first non-trivial exponential bounds on the size of the considered combinatorial structures. In a follow-up note, Tao and Sawin \cite{Blog2} showed that the bounds derived in \cite{ellenberg-gjswijt-2017} and \cite{naslund-sawin-2017} are exponentially optimal in the context of this polynomial method. Namely, no further exponential improvement can be obtained by more refined bounds on the slice rank of the adopted tensors, since the computed upper bounds coincide with the true values asymptotically to the first order in the exponent (that is, the bounds on the \emph{asymptotic} slice ranks, in the sense of \cite{jeroen-2018}, are tight). 

A problem which has a similar flavor, but possibly a different nature, is that of determining the exponential grow of \emph{trifferent} sets of ternary sequences. In this case one asks for the size of the largest subset of $\mathbb{F}_3^n$ with the property that any three distinct elements are simultaneously distinct in at least one coordinate. That is, they are projected onto $\mathbb{F}_3$ in at least one coordinate. This problem originates both in the context of information theory and in theoretical computer science, respectively as a problem of zero-error capacity under list decoding (or hypergraph capacity) or perfect hashing in a ternary alphabet. See \cite[Prob. 10.29]{csiszar-korner-book}, \cite{korner-simonyi}, \cite{dalai-iee} for further details.

If $T(n)$ is the size of a largest trifferent subset of $\mathbb{F}_3^n$, one can prove easily by induction that 
$$
T(n)\leq 2\left(\frac{3}{2}\right)^n\,.
$$
So, in this case there is already a simple non-trivial exponential upper bound on the size of the combinatorial structure, and it is rather natural to ask whether the slice rank method can be used to improve upon it.
Neglecting for a moment the details of how one might try to encode the trifference problem in the slice-rank method, one is easily led to the question of whether there exists at all tensors whose $n$-fold tensor powers have a slice rank which grows exponentially slower than $(3/2)^n$. 

In general, we ask what the smallest possible asymptotic slice rank  of a tensor can be, assuming the tensor is not a slice already.
The main result of this paper is to show that there is indeed a gap; any $k$-tensor in any field has either asymptotic slice rank 1 (i.e., it is a slice) or at least $k/(k-1)^{(k-1)/(k)}$. This can be interpreted as an extension of the fact that the standard asymptotic rank of matrices is either 1 or at least $2$.
In the case $k=3$, the found value is $3/2^{2/3}\approx 1.889>3/2$, so that no straight-forward application of the slice rank method can give improvements over known bounds for the trifference problem (but see Section \ref{Section3}).

In Section \ref{Section2} we prove our statement on the asymptotic slice rank of $k$-tensors. The main technical contribution is to show that all instances of  the method of Tao and Sawin (which depend on tensor representation) give a trivial lower bound if and only if the tensor has slice-rank 1.
In Section \ref{Section3}, taking again inspiration from the trifference problem, we add some comments on the limitations of our own result.
\subsection{Notation}
Following \cite{Blog2}, we consider finite-dimensional vector spaces, $V_1,\dots,V_k$, over a field $F$ and a basis $B_i=(b_{i,s})_{s\in S_i}$ for each $V_i$, $i\in[1,k]$, indexed by some finite set $S_i\subset \mathbb{Z}$. 
Given $\Gamma\subseteq S_1\times\dots \times S_k\subset \mathbb{Z}^k$, a $k$-tensor of $\bigotimes_{i=1}^k V_i$ will be defined as:
$$v=\sum_{(s_1,\dots,s_k)\in \Gamma} c_{s_1,\dots,s_k} b_{1,s_1}\otimes\dots\otimes b_{k,s_k}.$$
In case all the coefficients are nonzero, $\Gamma$ is said to be the support of $v$ with respect to the bases $B=\{B_1,\dots,B_k\}$.
For each $1\leq j\leq k$, we use the $j^{th}$ tensor product $\otimes_j:\ V_j \bigotimes_{i=1, i\not=j}^k V_i\rightarrow \bigotimes_{i=1}^k V_i$ as defined in \cite{Blog2} and let $\pi_j$ be the projection on the $j$-th coordinate. 

The notion of rank that will be discussed here is the following one:
\begin{Definition}
Tensors of the form $v_j\otimes_j v_{\hat{j}}$ for some $v_j\in V_j$ and $v_{\hat{j}}\in \bigotimes_{i=1,\ i\not=j}^k V_i$ have slice rank one and are said to be slice tensors.
The slice rank of an element of  $\bigotimes_{i=1}^k V_i$ is defined to be the least non negative integer $r$ such that $v$ is a linear combination of $r$ slice tensors.
\end{Definition}

In the note \cite{Blog2}, Terence Tao and William Sawin introduce a combinatorial way to study the slice rank of tensors. The key idea is to study the entropy of a set $\Gamma\subset \mathbb{Z}^k$ defined as follows:
$$H(\Gamma):=\sup_{(X_1,\dots,X_k)} \min(h(X_1),\dots,h(X_k)),$$
where $(X_1,\dots,X_k)$ ranges over the random variables taking values in $\Gamma$ and $h(X)$ is the Shannon entropy of the discrete variable $X$. In detail, 
$$h(X)=-\sum_{\alpha}p_{\alpha} \log(p_{\alpha})\,,$$
where $p_{\alpha}$ is the probability that $X=\alpha$ (we set $0\log 0 = 0$ by definition). Note that $h(X)=0$ if and only if $X$ is constant, so that $H(\Gamma)=0$ if and only if at least one of the coordinates is constant in $\Gamma$.

Using this notations, the upperbound on the slice-rank of tensors derived in \cite{Blog2} can be stated as follows.

\begin{Proposition}[\cite{Blog2}]\label{Upperbound}
Let $v$ be a $k$-tensor and let $\Gamma$ be its support with respect to the bases $B$. Then:
$$ \rank(v^{\otimes n})\leq \exp((H(\Gamma)+o(1))n).$$
\end{Proposition}
In \cite{Blog2}, the authors also provide a lower bound on the asymptotic slice rank of a $k$-tensor $v$. At this purpose, given total orderings $\sigma_1,\dots,\sigma_k$ for the finite sets $S_1,\dots,S_k$, consider the product ordering $\sigma=\sigma_1\times \sigma_2\times\dots \times \sigma_k$. Since $\sigma$ is a partial ordering, for any subset $\Gamma$ we can define $$\Gamma_{\sigma}:=\max_{\sigma}(\Gamma)$$
i.e. the set of maximal elements of $\Gamma$ with respect to $\sigma$.
The lowerbound of \cite{Blog2} can be states as follows:
\begin{Proposition}[\cite{Blog1}]\label{Lowerbound}
Let $v$ be a $k$-tensor and let $\Gamma$ be its support with respect to the bases $B$. Then, given a (product) ordering $\sigma$:
$$ \rank(v^{\otimes n})\geq \exp((H(\Gamma_{\sigma})+o(1))n).$$
\end{Proposition}
Proposition \ref{Lowerbound} lower bounds the slice rank of a tensor in terms of the entropy of antichain of maximal elements of its support with respect to some bases and ordering. For a given tensor, different bases and different orderings will in general give different lower bounds, and one might wonder what the best choice is. In particular, one might ask how small can the right hand side be, for some tensors, even for the optimal choice.
Our main result is that either the tensor is a slice or one can choose bases and ordering for which $H(\Gamma_{\sigma})\geq \log(k/(k-1)^{(k-1)/(k)})$. Using Proposition \ref{Lowerbound}, this leads to the following result.

 \begin{Theorem*}\label{thm1}
Let $v$ be a $k$-tensor that is not a slice. Then:
$$ \rank(v^{\otimes n})\geq (k/(k-1)^{(k-1)/(k)})^{n+o(n)}.$$
\end{Theorem*}
We prove this in two steps. First we show that, for any ordering $\sigma$ and finite set $\Gamma$, the quantity $H(\Gamma_{\sigma})$ is either zero or at least $\xi_k := \log(k/(k-1)^{(k-1)/(k)})$. The main task is then to show that, if $v$ is a $k$-tensor which is not a slice, then there exist bases $B$ and an ordering $\sigma$ with respect to which the tensor support $\Gamma$ satisfies $H(\Gamma_{\sigma})\not=0$. Hence, using Proposition \ref{Lowerbound}, $\rank(v^{\otimes n})\geq \exp(\xi_k(n+o(n)))$.

\section{Proof of Theorem \ref{thm1}}\label{Section2}

\begin{Lemma}\label{extSub}
Let $\Gamma\subseteq \mathbb{Z}^k$ with $k\geq 2$. Then $H(\Gamma)\not=0$ implies that there exists $\bar{\Gamma}\subseteq \Gamma$ such that $H(\bar{\Gamma})\not=0$ and $|\bar{\Gamma}|\leq k$. 
\end{Lemma}
\proof
We prove that, given $\Gamma$ such that $|\pi_1(\Gamma)|>1,\dots,|\pi_{k}(\Gamma)|>1$, there exists $\bar{\Gamma}\subseteq \Gamma$ such that $|\bar{\Gamma}|\leq k$ while still $|\pi_1(\bar{\Gamma})|>1,\dots,|\pi_{k}(\bar{\Gamma})|>1$, implying $H(\bar{\Gamma})>0$. We proceed by induction.

For the base case, $k=2$, let $(x_1,x_2)\in \Gamma$. Since $|\pi_1(\Gamma)|>1$ and $|\pi_2(\Gamma)|>1$, there exist $(y_1,y_2)$ and $(z_1,z_2)$ such that $y_1\not=x_1$ and $z_2\not=x_2$. If $y_2\not=x_2$ (or $z_1\not=x_1$) we can choose $\bar{\Gamma}$ to be $\{(x_1,x_2),(y_1,y_2)\}$ (reps.$\{(x_1,x_2),(z_1,z_2)\}$). Assuming $y_2=z_2$ and $z_1=x_1$, instead, we can choose $\bar{\Gamma}$ to be $\{(y_1,y_2),(z_1,z_2)\}$.

Assume now the statement if proved for $k-1$ and consider $\Gamma\subset \mathbb{Z}^k$. Let $\tilde{\Gamma}\subset \Gamma$ be the smallest subset of $\Gamma$ such that $|\pi_1(\tilde{\Gamma})|>1,\dots,|\pi_{k-1}(\tilde{\Gamma})|>1$. Because of the inductive hypothesis (applied to the projection of $\tilde{\Gamma}$ on the first $k-1$ coordinates), we have that also $|\tilde{\Gamma}|\leq k-1$. If also $|\pi_k(\tilde{\Gamma})|>1$, then $\tilde{\Gamma}$ is already a subset of $\Gamma$ that satisfies the required properties. Otherwise, since $\pi_k(\Gamma)>1$, there exists $x\in \Gamma$ so that, set $\bar{\Gamma}:=\{x\}\cup \tilde{\Gamma}$, we have  $|\pi_1(\bar{\Gamma})|>1,\dots,|\pi_{k}(\bar{\Gamma})|>1$. Since $|\tilde{\Gamma}|\leq k-1$ we also have that $|\bar{\Gamma}|\leq k$.
\endproof
\begin{Proposition}\label{Gap1}
Let $\Gamma$ be a finite subset of $\mathbb{Z}^k$. Then denoted by $\xi_k:=\log(k/(k-1)^{(k-1)/k})$, for any ordering $\sigma$, $H(\Gamma_{\sigma})\not\in\  (0,\xi_k)$. 
Moreover, there exist $\Gamma\subseteq \mathbb{Z}^k$ and $\sigma$ such that $H(\Gamma_{\sigma})=\xi_k$. \end{Proposition}
\proof
Let us consider an ordering $\sigma$ such that $H(\Gamma_{\sigma})\not=0$. Because of Lemma \ref{extSub}, we may suppose $|\Gamma_{\sigma}|\leq k$.
It follows from the definition that $H(\Gamma_{\sigma})\geq \min(h(X_1),\dots, h(X_k))$ where $(X_1,\dots,X_k)$ is the uniformly distributed random variable on $\Gamma_{\sigma}$. The rest is just an application of the data processing inequality for the entropy; we write the details for readers with a different background.
For any $i\in [1,k]$, denote by $p_{\alpha}$ the probability of the event $X_i=\alpha$, so that $h(X_i)=-\sum_{\alpha\in S_i}p_{\alpha} \log(p_{\alpha})$.
Since $|\pi_i(\Gamma_\sigma)|\not=1$, there exists $\bar{\alpha}\in S_i$ such that $1/k\leq p_{\bar{\alpha}}\leq 1/2$.
Then, it follows from Jensen's inequality after simple algebraic manipulations that
$$-\sum_{\alpha\in S_i}p_{\alpha} \log(p_{\alpha})\geq -p_{\bar{\alpha}}\log(p_{\bar{\alpha}})-\left(\sum_{\alpha\in S_i, \alpha\not=\bar{\alpha}}p_{\alpha}\right) \left(\log(\sum_{\alpha\in S_i, \alpha\not=\bar{\alpha}}p_{\alpha})\right).$$
Since $1/2\geq p_{\bar{\alpha}}\geq 1/k$ and the function $-x\log(x)-(1-x)\log(1-x)$ is monotonic in $[0,1/2]$, we have that
$$h(X_i)\geq -1/k\log(1/k)-(k-1)/k\log((k-1)/k)=\log(k/(k-1)^{(k-1)/k}).$$
Summing up we obtain that
$$H(\Gamma_{\sigma})\geq \min(h(X_1),\dots, h(X_k))\geq \log(k/(k-1)^{(k-1)/k}).$$
On the other hand, it is easy to check that the set of $k$ points in $ \mathbb{Z}^k$
$$\Gamma:=\{(2,1\dots,1),(1,2,1,\dots,1),\dots,(1,1,\dots,2)\}$$
is such that $H(\Gamma_{\sigma})=H(\Gamma)=\log(k/(k-1)^{(k-1)/k})$, where $\sigma$ is the usual product ordering ($2>1$ in any coordinate).
\endproof
It follows from Proposition \ref{Gap1} and Proposition \ref{Lowerbound} that:
\begin{Corollary}\label{Cor2}
Let $v$ be a $k$-tensor and let $\Gamma$ be its support with respect to the bases $B$. If there exists an ordering $\sigma$ such that $H(\Gamma_{\sigma})\not=0$, we have that:
$$ \rank(v^{\otimes n})\geq \exp(\xi_k(n+o(n)).$$
\end{Corollary}
Due to Corollary \ref{Cor2}, we would like to characterize the tensors $v$ whose support $\Gamma$ with respect to any bases and ordering satisfies $H(\Gamma_{\sigma})=0$, and show that they must be slice. 
Toward that goal, we first investigate the entropy of \emph{sections} of a set $\Gamma$, defined as follows. Given $\Gamma\subseteq\mathbb{Z}^k$, $I=\{i_1,i_2,\ldots,i_t\}\subseteq[1,k]$ and $x=(x_1,x_2,\ldots,x_t)$, let $M_I^x\subseteq \Gamma$ be the subset of elements with $j$-th component $x_{i_j}$. Let then $\Gamma_I^x\subseteq \mathbb{Z}^{k-t}$ be the projection of $M_I^x$ on the coordinates $[1,k]\setminus I$.
\begin{Lemma}\label{ProdCartesiano}
Let $\Gamma$ be a finite subset of $\mathbb{Z}^k$ and let $P_i=\pi_i(\Gamma)$. Then $H((\Gamma_I^x)_{\sigma})=0$ for any $I\subseteq [1,k]$, any $x\in \bigtimes_{i\in I}P_i$ and any ordering $\sigma$ if and only if $\Gamma=\bigtimes_{i\in[1,k]} P_i$.
\end{Lemma}
\proof
If $\Gamma=\bigtimes_{i\in[1,k]} P_i$, then for any $I\subseteq [1,k]$ and $x\in \bigtimes_{i\in I}P_i$ we have that $\Gamma_I^x$ is the cartesian product of $P_i $ such that $i\not\in I$. Therefore, $\Gamma_I^x$ has a maximum with respect to any ordering $\sigma$, which implies that $H((\Gamma_I^x)_{\sigma})=0$.

Let us suppose $\Gamma$ is not the cartesian product of the sets $P_i$. Up to permutation of the coordinates, this means that there exist $x=(x_1,\dots,x_k)$ and $z=(z_1,x_2,\dots,x_k)$ such that $z_1$ and $x_1$ are distinct elements of $P_1$, $x\in \Gamma$ and $z\not \in \Gamma$.  
Let us now consider in $\Gamma$ an element $\bar{z}$ with $\pi_1(\bar{z})=z_1$ that differs from $x$ in the minimum number, say $t$, of coordinates; since $\bar{z}\neq z\not\in \Gamma$, $t\geq 2$. We can assume, up to a rearrangement of the coordinates, that $\bar{z}=(z_1,z_2,\dots, z_t,x_{t+1},\ldots,x_{k})$, with $z_i\neq x_i$ for $i\in[1,t]$. Set $x'=(x_{t+1},\dots,x_{k})$ and $I=[t+1,\dots,k]$; on the section $\Gamma_I^{x'}$, we consider the product ordering $\sigma$ on $\mathbb{Z}^t$ such that in the first coordinate $z_{1}$ is the largest element and $x_{1}$ the second largest one, let us write $z_{1}>x_{1}>\dots$ , while in the other coordinates, using the same notation\footnote{In the rest of the paper we always use this notation to specify the largest and the second largest elements of the ordering.\label{footnoteordering}}, $x_{2}>z_{2}>\cdots,\cdots, x_t>z_t>\cdots$.

Note that $(z_{1},\dots,z_{t})$ belongs to $(\Gamma_I^{x'})_{\sigma}$ because no element can majorize it, since it would need to have some $i$-th coordinate equal to $x_i$ while $\bar{z}$ differs from $x$ in the minimum number of coordinates. But $(x_{1},\dots,x_{t})$ also belongs to $(\Gamma_I^{x'})_{\sigma}$. Since $z_i\neq x_i$ for $i\in [1,t]$, we deduce that $H((\Gamma_I^{x'})_{\sigma})\not=0$
\endproof

\begin{Lemma}\label{TrasportoAnticatena}
Let $\Gamma$ be a finite subset of $\mathbb{Z}^{k}$, let $P_i=\pi_i(\Gamma)$.
Suppose that the following conditions hold:
\begin{itemize}
\item set $I=[1,k-d]$, there exists $x=(x_1,\dots,x_{k-d})\in P_1\times \dots \times P_{k-d}$ and an ordering $\sigma$ with $H((\Gamma_I^{x})_{\sigma})\not=0$;
\item for any $I'\subset I$ and for any ordering $\alpha$, $H((\Gamma_{I'}^{x'})_{\alpha})=0$ where $x'$ is the restriction of $x$ on $I'$.
\end{itemize} 
Then set $\bar{x}=(x_2,\dots,x_{k-d})$ the restriction of $x$ on $\bar{I}=[2,k-d]$, we have that $(\Gamma_I^{x})_{\sigma}\subseteq \Gamma_{I}^{(z_1,\bar{x})}$ for any $z_1\in P_1$.
\end{Lemma}
\begin{figure}
\centering
\includegraphics[scale=0.7]{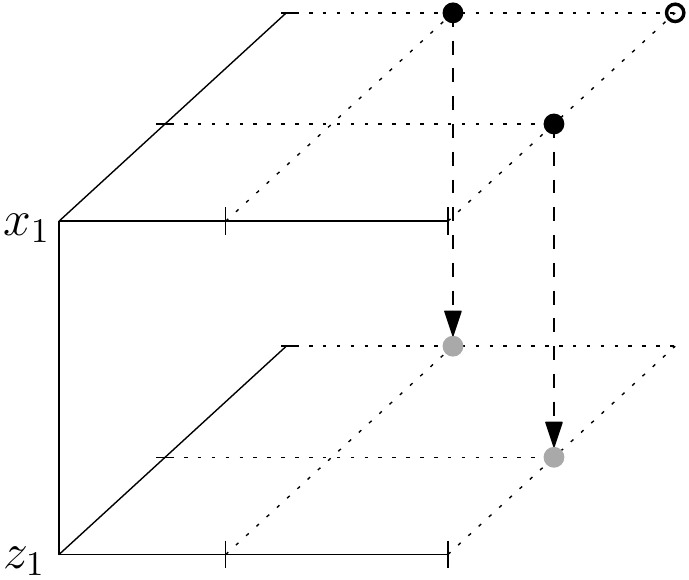}
\caption{Graphical representation of Lemma \ref{TrasportoAnticatena} for the case where $d=2$. Coordinates in $\bar{I}=[2,k-d]$, which are fixed to $\bar{x}$, are omitted.}
\label{fig1}
\end{figure}
\proof
If $|P_1|=1$ the statement is trivial, so let $z_1\in P_1$ be different from $x_1$.

First of all we prove that $\Gamma_{I}^{(z_1,\bar{x})}$ is not empty. Let us consider an element $z \in P_1\times \dots \times P_{k-d}$ with $\pi_1(z)=z_1$ and non empty $\Gamma_I^{z}$ which differs from $x$ in the minimum number of coordinates. We would like to prove that $z=(z_1,\bar{x})$. We can assume, up to a rearrangement of the coordinates, that $z=(z_1,z_2,\dots,z_t,x_{t+1},\dots,x_{k-d})$ where $x_i\not=z_i$ for any $i\in [1,t]$. We set $x'=(x_{t+1},\dots,x_{k-d})$ and $I'=[t+1,d]$. Let $\alpha$ be the ordering on $\Gamma_{I'}^{x'}$ that coincide with $\sigma$ on $P_{k-d+1}\times \dots \times P_k$ and such that (see footnote \ref{footnoteordering}) $z_1>x_1>\cdots$, $x_2>z_2>\cdots,\dots, x_{t}>z_{t}>\cdots$. We have that $\{(z_1,\dots,z_t)\}\times (\Gamma_I^{z})_{\sigma}\subseteq (\Gamma_{I'}^{x'})_{\alpha}$ because no element can majorize it, since it would need to have some $i$-th coordinate equal to $x_i$ while $z$ differs from $x$ in the minimum number of coordinates. In the case $t\geq 2$, $\Gamma_{I}^{(z_1,\bar{x})}$ would be empty and hence we also have that $\{(x_1,\dots,x_t)\}\times (\Gamma_I^{x})_{\sigma}\subseteq (\Gamma_{I'}^{x'})_{\alpha}$. 
Since $H((\Gamma_I^{x})_{\sigma})\not=0$, for $t\geq 2$ none of the last $d$ coordinates of points in $(\Gamma_{I'}^{x'})_{\alpha}$ is constant. Also none of the first $t$ coordinates would be constant because $x_i\not=z_i$ for any $i\in [1,t]$ and this would imply $H((\Gamma_{I'}^{x'})_{\alpha})\not=0$. Since instead, $H((\Gamma_{I'}^{x'})_{\alpha})=0$ we necessarily have that $t=1$, $x'=\bar{x}=(x_2,\dots,x_{k-d})$, $z=(z_1,\bar{x})$ and $\Gamma_I^{(z_1,\bar{x})}\not=\emptyset$.

Let now consider $\beta_1$ to be an ordering on $\mathbb{Z}$ such that $x_1>z_1>\cdots$. 
Then consider over $\mathbb{Z}^{d+1}$ the product ordering $\beta=\beta_1\times \sigma$. Note that any element in $\{x_1\}\times (\Gamma_I^{x})_{\sigma}$ is also in $(\Gamma_{\bar{I}}^{\bar{x}})_{\beta}$, since it cannot be majorized by any other one. So, the assumptions $H((\Gamma_I^{x})_{\sigma})\neq 0$ and  $H(\Gamma_{\beta})=0$ imply that indeed $(\Gamma^{\bar{x}}_{\bar{I}})_{\beta}=\{x_1\}\times (\Gamma_I^{x})_{\sigma}$.
We now swap $x_1$ and $z_1$; let $\gamma_1$ be an ordering on $\mathbb{Z}$ such that $z_1>x_1>\cdots$, and set similarly $\gamma=\gamma_1\times \sigma$. Now, for any $y\in (\Gamma_I^{x})_{\sigma}$, we note that $(x_1, y)\notin (\Gamma_{\bar{I}}^{\bar{x}})_{\gamma}$ if and only if $(z_1,y)\in  (\Gamma_{\bar{I}}^{\bar{x}})_{\gamma}$, for no other element can majorize it. So, either $(z_1,y)\in (\Gamma_{\bar{I}}^{\bar{x}})_{\gamma}$ or $(x_1,y)\in (\Gamma_{\bar{I}}^{\bar{x}})_{\gamma}$. 
Since $H((\Gamma_I^{x})_{\sigma})\neq 0$, none of the last 
$d$ coordinates of points in $(\Gamma_{\bar{I}}^{\bar{x}})_{\gamma}$ is constant. Hence, $H((\Gamma_{\bar{I}}^{\bar{x}})_{\gamma})=0$ implies that the first coordinate is constant in  $(\Gamma_{\bar{I}}^{\bar{x}})_{\gamma}$. But, since $\Gamma_I^{(z_1,\bar{x})}$ is not empty, $(\Gamma_{\bar{I}}^{\bar{x}})_{\gamma}$ must contain some points in $\{z_1\}\times \Gamma_I^{(z_1,\bar{x})}$, in particular all points in $\{z_1\}\times (\Gamma_I^{(z_1,\bar{x})})_{\sigma}$. So, $(\Gamma_{\bar{I}}^{\bar{x}})_{\gamma}\subseteq\{x_1\}\times \Gamma_I^{(z_1,\bar{x})}$ and any $y\in (\Gamma_{I}^{x})_{\sigma}$ is also contained in $\Gamma_{I}^{(z_1,\bar{x})}$.
\endproof

\begin{Lemma}\label{CreareBuco}
Let $v$ be a $k$-tensor and let $\Gamma_B$ be its support with respect to the bases $B$. Let us suppose that there exist $x_i,z_i$ in $\pi_i(\Gamma_B)$ and $y\in (\Gamma_B)_i^{x_i}\cap (\Gamma_B)_i^{z_i}$. Then there exist bases $B'$ such that:
\begin{itemize}
\item $|B'|=|B|$; 
\item$(\Gamma_B)_i^{x_i}=(\Gamma_{B'})_i^{x_i}$;
\item $y\not\in (\Gamma_{B'})_i^{z_i}$.
\end{itemize}
\end{Lemma}
\proof
We assume, without loss of generality, that $i=1$ and we set $y=(y_2,\dots,y_k)$. As usual we have the following expression for $v$:
$$v=\sum_{(s_1,\dots,s_k)\in \Gamma_{B}} c_{s_1,s_2\dots,s_k} b_{1,s_1}\otimes \dots \otimes b_{k,s_k}$$
where all coefficients are non-zero.

Now we would like to proceed with the Gaussian elimination respect to the basis $B_1=(b_{1,j}| j\in S_j)$. We consider the vectors $b'_{1,x_1}=c_{x_1,y_2\dots,y_k} b_{1,x_1}+c_{z_1,y_2\dots,y_k} b_{1,z_1}$ and $b'_{1,z_1}=b_{1,z_1}$ and we set $B'=(B\setminus \{b_{1,x_1},b_{1,z_1}\})\cup \{b'_{1,x_1},b'_{1,z_1}\}$. Then we have the following expression for $v$:
$$v=\sum_{(s_1,\dots,s_k)\in \Gamma_{B'}} c'_{s_1,s_2\dots,s_k} b_{1,s_1}\otimes \dots \otimes b_{k,s_k}$$
where all coefficients are non-zero. We have that:
$$c'_{x_1,s_2\dots,s_k}=c_{x_1,s_2\dots,s_k}/c_{x_1,y_2\dots,y_k}.$$
Therefore $(\Gamma_{B'})_1^{x_1}=(\Gamma_{B})_1^{x_1}$.
Similarly we have that:
$$c'_{z_1,s_2\dots,s_k}=c_{z_1,s_2\dots,s_k}+c_{x_1,s_2\dots,s_k}(-c_{z_1,y_2\dots,y_k}/c_{x_1,y_2\dots,y_k})$$
and hence $y\not\in (\Gamma_{B'})_1^{z_1}$.
Therefore we have found bases $B'$ such that $|B'|=|B|$, $(\Gamma_{B'})_1^{x_1}=(\Gamma_{B})_1^{x_1}$ and $y\not\in (\Gamma_{B'})_1^{z_1}$.
\endproof
As a consequence of the previous lemmas we can now prove the following proposition.
\begin{Proposition}\label{NSliceHN0}
Let $v$ be a $k$-tensor that is not a slice. Then there exist bases $B$ and an ordering $\sigma$ such that, denoted by $\Gamma$ the support of $v$ respect to $B$, we have that $H(\Gamma_{\sigma})\not=0$.
\end{Proposition}
\proof
We first note that we can choose among the bases $B$ with minimal $|B|$ one for which there exists a section with non-zero entropy. In fact, because of Lemma \ref{CreareBuco} we can choose a basis for the which the support $\Gamma$ is not the cartesian product of the sets $P_i=\pi_i(\Gamma)$. Therefore, according to Lemma \ref{ProdCartesiano}, we can assume that there exist a section, say $\Gamma_I^x$, and an ordering $\sigma$ with $H((\Gamma_I^x)_{\sigma})\not=0$. 
Among such bases with minimum $|B|$, we choose bases $B$, an $I$ and an $x$ so that the dimension $d$ of $\Gamma_I^x$ is the maximal one. 

Let us suppose, by contradiction, that $d<k$. Up to a permutation of the coordinates we can assume $I=[1,k-d]$ and write $x=(x_1,\dots,x_{k-d})$. Since $v$ is not a slice, there exist $z_{1}$ distinct from $x_{1}$ in $P_{1}$. Set $z=(z_1,x_2,\dots,x_{k-d})$. 
Let also $\bar{I}=[2,k-d]$ and $\bar{x}=(x_2,\dots,x_{k-d})$. Applying then Lemma \ref{TrasportoAnticatena} to the section $\Gamma_{\bar{I}}^{\bar{x}}$ we deduce that $(\Gamma_I^x)_{\sigma}\subseteq \Gamma_I^{(z_1,\bar{x})}=\Gamma_I^z$. Therefore, given $y\in (\Gamma_I^x)_{\sigma}$, we 
 have that $y \in  \Gamma_{I}^{x}\cap \Gamma_{I}^{z}$ or, equivalently, $(\bar{x},y)\in \Gamma_{1}^{x_{1}}\cap \Gamma_{1}^{z_{1}}$. 
Therefore, as a consequence of Lemma \ref{CreareBuco}, there exists bases $B'$ such that, denoted by $\Gamma_{B'}$ the support of $v$ respect to $B'$:
\begin{itemize}
\item[1)] $|B'|=|B|$; 
\item[2)] $\Gamma_{1}^{x_{1}}=(\Gamma_{B'})_{1}^{x_{1}}$ and hence $\Gamma_I^x=(\Gamma_{B'})_I^x$;
\item[3)] $(\bar{x},y)\not\in (\Gamma_{B'})_{1}^{z_{1}}$ that is $y\not\in (\Gamma_{B'})_{I}^{z}$.
\end{itemize}
Because of $2)$, if we consider the ordering $\sigma$ defined above, we still have that $H(((\Gamma_{B'})_{I}^x)_{\sigma})\not=0$ and, because of the maximality of $d$, we have that $H(((\Gamma_{B'})_{I'}^{x'})_{\alpha})=0$ for any $I'\subset I$ and for any ordering $\alpha$, where $x'$ is the restriction of $x$ on $I'$. Due to the minimality of $|B|=|B'|$, $z_{1}\in \pi_{1}(\Gamma_{B'})$ otherwise, using the notation of Lemma \ref{CreareBuco}, we could omit $b_{1,z_{1}}$ from $B'$ obtaining a smaller bases.
Now we can apply again Lemma \ref{TrasportoAnticatena} to the section $(\Gamma_{B'})_{\bar{I}}^{\bar{x}}$ obtaining that $((\Gamma_{B'})_{I}^x)_{\sigma}\subseteq (\Gamma_{B'})_I^z$. But this is in contradiction with $3)$ because $y$ is in $(\Gamma_{I}^x)_{\sigma}=((\Gamma_{B'})_{I}^x)_{\sigma}$ but not in $(\Gamma_{B'})_{I}^{z}$. It follows that $d=k$ and hence there exist bases $B$ and an ordering $\sigma$ such that, denoted by $\Gamma$ the support of $v$ respect to $B$, we have that $H(\Gamma_{\sigma})\not=0$.
\endproof

Combined with Corollary \ref{Cor2}, Proposition \ref{NSliceHN0} proves Theorem \ref{thm1}. 

\section{On variations of the method}\label{Section3}
As mentioned in the introduction, our study was initially motivated by the trifference problem. An immediate consequence of our result is that one cannot hope to derive bounds smaller than $1.889^n$ on the size of trifferent sets by applying the slice rank method \emph{in a straight-forward way}, that is using a $3$-tensor which is a tensor power and whose coordinates are indexed by elements of $\mathbb{F}_3^n$ as done for the capset problem. 

However, this does not imply that the slice rank method cannot be used at all by means of more elaborate applications. We show that, for example, one can actually prove a bound of $3^{n/2}\approx 1.732^n$ on the size of trifferent sets using the polynomial method with just $2$-tensors whose rows and columns are indexed by pairs of distinct sequences (an instance of the slice rank method which boils down to the original method of Haemers \cite{haemers-1978} for bounding the graph capacity).
Of course this is way worse than the best known bound of $2(3/2)^n$ mentioned in the introduction, but is suffices to show that the gap proved for $k$-tensors should \emph{never} be interpreted to mean that no use can be made in general of the slice rank method for a given problem.

Let $A\subset \{1,\omega,\omega^2\}^n$ be a trifferent set, where $\omega=e^{i2\pi/3}$. Let for simplicity $A^{(2)}$ be the set of $|A|(|A|-1)/2$ unordered pairs of distinct elements of $A$.

For $(x,y)\in A^{(2)}$, consider the function $f_{x,y}:A^{(2)}\to \mathbb{C}$ defined by
$$
f_{x,y}(z,t)=\prod_{i=1}^n (x_i+y_i+z_i)(x_i+y_i+t_i)\,.
$$
If $(x,y)=(z,t)$, then 
\begin{align*}
f_{x,y}(z,t) & =\prod_{i=1}^n (2x_i+y_i)(x_i+2y_i) \\
& \neq 0.
\end{align*}
If $(x,y)\neq (z,t)$ then either $(x,y,z)$ or $(x,y,t)$ is a trifferent triplet and hence either $(x_i+y_i+z_i)=0$ for some $i$ or $(x_i+y_i+t_i)=0$ for some $i$. So
$$
f_{x,y}(z,t)=0\,,\quad (z,t)\neq(x,y)\,.
$$
This implies that the functions $f_{x,y}$ with $(x,y)\in A^{(2)}$ are linearly independent, because if
$$
\sum_{x,y}a_{x,y}f_{x,y} = 0
$$
then computing the left hand side on $(z,t)$ we find $a_{z,t}=0$.

But we can write
$$
f_{x,y}(z,t)=\prod_{i=1}^n ((x_i+y_i)^2+(x_i+y_i)(z_i+t_i)+z_it_i)\,,
$$
which can be expanded as the sum of $3^n$ terms of the form 
$$
c\prod_{i_1}^n(z_i+t_i)^{\alpha_i}(z_it_i)^{\beta_i}
$$
with $\alpha_i,\beta_i\in \{0,1\}$, $\alpha_1+\beta_i\leq 1$.
So, the functions $f_{x,y}$ live in a space of dimension at most $3^n$. This implies that asymptotically
\begin{align*}
|A| & \leq (\sqrt{3}+o(1))^n\\
& \approx (1.7321)^n.
\end{align*}

Note that the above procedure can be interpreted as an instance of the standard polynomial method but also as an instance of the slice rank method which, when applied to $2$-tensors, is essentially equivalent to the original method of Haemers \cite{haemers-1978}.

\section*{Acknowledgements}

This research was partially supported by Italian Ministry of Education under Grant PRIN 2015 D72F16000790001. Helpful discussions with Jaikumar Radhakrishnan and Venkatesan Guruswami are gratefully acknowledged.

\bibliographystyle{plainurl}
\bibliography{slice-rank}

%\begin{thebibliography}{Z}
%\bibitem{croot-lev-pach-2016}

%\bibitem{Blog1}Terence Tao, William Sawin. Notes on the slice rank of tensors, 2016. URL: https://terrytao.wordpress.com/2016/08/24/notes-on-the-slice-rank-of-tensors/ 
%\bibitem{Blog2}Terence Tao, A symmetric formulation of the Croot--Lev--Pach--Ellenberg-Gijswijt capset bound, 2016. URL:https://terrytao.wordpress.com/2016/05/18/a-symmetric-formulation-of-the-croot-lev-pach-ellenberg-gijswijt-capsetbound
%\end{thebibliography}
\end{document}